\documentclass[10pt]{article}
\usepackage{amssymb,amsmath,amsthm}
\newcommand{\ie}{\emph{i.e.}}
\newcommand{\eg}{\emph{e.g.}}
\newcommand{\cf}{\emph{cf}}
\newcommand{\ad}{\emph{ad}}
\newcommand{\etal}{\emph{et al.}}
\newcommand{\Com}{\mathbb{C}}
\newcommand{\Real}{\mathbb{R}}
\newcommand{\Nat}{\mathbb{N}}
\newcommand{\s}{L}
\newcommand{\sloc}{L_\mathrm{loc}}
\newcommand{\Smooth}{C}
\newcommand{\Dom}{D}
\newcommand{\id}{1}
\newcommand{\Span}{\mathop{\mathrm{span}}}
\newcommand{\Hilbert}{\mathcal{H}}
\newcommand{\esssup}{\mathop{\mathrm{ess\;\!sup}}}
\newcommand{\essinf}{\mathop{\mathrm{ess\;\!inf}}}

\newcommand{\demi}{\frac{1}{2}}
\newcommand{\re}{\mathop{\mathrm{Re}}}
\newcommand{\op}{S}
\newcommand{\opf}{s}
\newcommand{\ope}{\mathfrak{S}}
\newcommand{\Hbasic}{\langle\mathrm{H1}\rangle}
\newcommand{\Hbasicbis}{\langle\mathrm{H1'}\rangle}
\newcommand{\Hsign}{\langle\mathrm{H2}\rangle}
\newenvironment{Assumption}[1]
{\begin{description}\item[$#1$]\it\ }
{\end{description}}
\newtheorem{Lemma}{Lemma}
\newtheorem{Proposition}{Proposition}
\newtheorem{Corollary}{Corollary}
\newtheorem{Theorem}{Theorem}
\theoremstyle{definition}
\newtheorem*{Remark}{Remark}
\newtheorem*{Remarks}{Remarks}
\begin{document}
%
\title{{\Large\textbf{
Instability results for the damped wave equation in unbounded domains}}
}
\author{
P.~Freitas 
\ and \ 
D.~Krej\v{c}i\v{r}\'{\i}k%
\footnote{
On leave from
\emph{Department of Theoretical Physics,
Nuclear Physics Institute, Academy of Sciences,
250\,68 \v{R}e\v{z} near Prague, Czech Republic}.
}
}
\date{\footnotesize
\begin{quote}
\begin{center}
\emph{
Departamento de Matem\'atica, Instituto Superior T\'ecnico, \\
Av. Rovisco Pais, 1049-001 Lisboa, Portugal 
\smallskip \\ 
\emph{E-mail:}
pfreitas@math.ist.utl.pt, 
dkrej@math.ist.utl.pt
}
\end{center}
\end{quote}
June 2, 2004
}
\maketitle
%
%
\begin{abstract}
\noindent 
We extend some previous results for the damped wave equation in 
bound\-ed domains in $\Real^{d}$ 
to the unbounded case. In particular, we 
show that if the damping term is of the form $\alpha a$ with 
bounded $a$ taking on negative values on a set of positive 
measure, then there will always exist unbounded solutions for 
sufficiently large positive~$\alpha$.

In order to prove these results, we 
generalize some existing results on the asymptotic behaviour of 
eigencurves of one-parameter families of Schr\"{o}dinger operators to the
unbounded case, which we believe to be of interest in their own right. 
\end{abstract}
%
\newpage
%
%
\section{Introduction}
%
Let~$\Omega$ be a  
domain in $\Real^d$ (bounded or unbounded), with $d \geq 1$,
and let~$\op_0$ denote the self-adjoint operator generated on~$\s^2(\Omega)$
by the differential expression $\ope_0 := -\Delta + b$,
subject to Dirichlet boundary conditions,
where~$b$ is a real-valued measurable function on~$\Omega$;
$\opf_0$ will denote the associated quadratic form.   
In this paper, we are interested in the instability of the solutions
$
  t \mapsto \psi(t,\cdot) 
  \in \Smooth^0\big([0,\infty);\Dom(\op_0)\big)
  \cap \Smooth^1\big([0,\infty);\Dom(\opf_0)\big)
  \cap \Smooth^2\big([0,\infty);\s^2(\Omega)\big) 
$  
of the following initial value problem 
\begin{equation}\label{wave.eq}
\left\{ 
\begin{aligned}
  \psi_{tt} + \alpha \, a \, \psi_{t} + \op_0 \psi & = 0 
  \,, 
  \\
  \psi(0,\cdot) & = \phi_1 \in \Dom(\op_0)
  \,,
  \\
  \psi_{t}(0,\cdot) & = \phi_2 \in \Dom(\opf_0)
  \,,
\end{aligned}
\right.
\end{equation}
where~$\alpha \geq 0$ is a parameter (not necessarilly small)
and~$a$ is a real-valued measurable function on~$\Omega$.

In the case where the damping $\alpha a$ remains non-negative and bounded, the 
asymptotic behaviour of solutions of~(\ref{wave.eq}) is well
understood~\cite{Zuazua}. However, there are situations for which this
will not be the case and the damping term will change sign, precluding 
the usage of energy methods and other tools which rely on the positivity of
the functional
\[
\int_{\Omega} a(x) \, |u(x)|^{2} \, dx \,.
\]
Such a situation might arise, for instance, when linearizing semilinear
damped wave equations around a stationary solution -- see, for
instance~\cite{DFG}. Because of the change in sign, it is not clear 
what will now happen to the stability of these solutions.

In 1991 Chen \etal~\cite{CFNS}
conjectured that for bounded intervals and under certain
extra conditions on the damping the trivial solution of~(\ref{wave.eq}) would
remain stable. This was
disproved in 1996 by the first author who showed that in the case of
bounded $\Omega$ in $\Real^{d}$
this sign-changing condition is sufficient to cause the existence of
unbounded solutions of~(\ref{wave.eq}), 
provided that the $L^{\infty}$-norm of the damping 
is large enough~\cite{Freitas_1996}. Heuristically, this behaviour can
be understood from the fact that, when $a$ changes sign, as the parameter
$\alpha$ increases equation~(\ref{wave.eq}) (formally) 
approaches a backward--forward heat equation and thus one does expect the
appearance of points in the spectrum on the positive side of the real axis.
On the other hand, and still at the heuristic level, note that while 
for bounded domains the result is not unexpected from the point 
of view of geometric optic rays either, for unbounded domains this is 
not as clear. Think, for instance, of the case of $\Omega = \Real^{d}$ and assume 
that $a(x)$ is negative inside the unit ball, positive outside, and goes to 
zero as $|x|$ goes to infinity. 
Then Theorem~\ref{Thm.instability} below still gives 
that for large enough $\alpha$ the trivial solution 
of~(\ref{wave.eq}) becomes unstable.

The purpose of the present paper is to extend the results for the bounded 
case to unbounded domains and, in 
fact, the idea behind the results given here is the same as that used
in~\cite{Freitas_1996}. The key point
is the observation (\cf~Lemma~\ref{Lem.basic})
that it is possible to obtain information
on the real part of the spectrum of the damped wave equation operator  
by studying some spectral properties of the one-parameter family of self-adjoint 
Schr\"odinger operator~$\op_\mu$ generated on~$\s^2(\Omega)$
by the differential expression 
\begin{equation}\label{diff.expression} 
  \ope_\mu := -\Delta + b + \mu \, a 
  \,, 
\end{equation}
where~$\mu$ is a real parameter. However, the situation now becomes
more complex due to the presence of essential spectrum, and, in 
particular, requires the extension of some results for~$\op_\mu$ to 
this setting which we believe to be interesting in their own right -- 
see Appendix~\ref{App.Asymptotes}.
 
To state the main instability result,
we now need to introduce some notation and basic assumptions.  
Throughout the paper, we assume that the damping coefficient~$a$ 
is bounded and that the potential~$b$ is bounded from below,
does not have local singularities but can grow at infinity, namely,
\begin{Assumption}{\Hbasic}   
$a\in\s^\infty(\Omega)$, \ $b\in\sloc^\infty(\Omega)$ 
\ and \ $b_\mathrm{min} := \essinf b > -\infty$\,.
\end{Assumption}
We also use the notation 
$$
  a_\mathrm{min} := \essinf a
  \qquad\mbox{and}\qquad 
  a_\mathrm{max} := \esssup a
  \,.
$$
Under the stated assumptions, the initial value problem~(\ref{wave.eq})
has a unique solution (\cf~Corollary~\ref{Corol.solution.wave} 
in Appendix~\ref{App.C0}).
The main result of the paper implies that 
if the damping~$a$ is sign-changing then
there are initial conditions of~(\ref{wave.eq}) 
for which the corresponding solution are unbounded
for sufficiently large~$\alpha$. 
\begin{Theorem}\label{Thm.instability}
Assume~$\Hbasic$. 
If~ $a_\mathrm{min}<0$, then there exists~$\alpha_0>0$
such that the system~\emph{(\ref{wave.eq})} is unstable 
for all~$\alpha>\alpha_0$. 
\end{Theorem} 

The instability property is obtained by studying a spectral problem
for a non-self-adjoint operator~$A_\alpha$  
generated by the matrix-differential expression
\begin{equation}\label{wave.op.formal}
  \begin{pmatrix}
    0       & 1 \\
    -\ope_0 & -\alpha\,a 
  \end{pmatrix} ,
\end{equation}
which appears in a first-order evolution equation
(\cf~(\ref{wave.eq.bis})) equivalent to~(\ref{wave.eq}),
and proving that under these conditions
there are positive points in its spectrum (Theorem~\ref{Thm.main}). 
If the point in the spectrum is an eigenvalue, then this result shows
that there are initial conditions of~(\ref{wave.eq}) 
for which the corresponding solution grows exponentially.
If the point belongs to the essential spectrum,
then Theorem~\ref{Thm.instability} follows as a consequence 
of the global instability property (GI-property) 
for a local semi-dynamical system
associated with~(\ref{wave.eq}), \cf~\cite{Sola-Morales}.
In order to apply the last result,
we need to ensure that~$A_\alpha$ generates a $C_0$-semigroup
(\cf~Appendix~\ref{App.C0}).

The organization of the paper is as follows.
In the preliminary Section~\ref{Sec.Preliminary},
we define properly the operators~$\op_\mu$ and~$A_\alpha$,
and state some basic properties.
The real spectrum of~$A_\alpha$ is investigated
in Section~\ref{Sec.Real}, where we also prove Theorem~\ref{Thm.main}
which implies the instability result of Theorem~\ref{Thm.instability}.
Section~\ref{Sec.Ess} is devoted to a more precise
discussion of the real essential spectrum of~$A_\alpha$ 
in unbounded domains if the damping~$a$ is asymptotically constant;
we also discuss there the relation between the essential spectrum
and the form of~$\Omega$ or the behaviour of~$b$. In 
Section~\ref{remarks}, we discuss the application of the results 
obtained to the semilinear case, and consider some open questions.
The paper is concluded by two appendices,
where we establish the asymptotics of eigenvalues of~$\op_\mu$ 
below the essential spectrum for large~$|\mu|$ 
and prove the $C_0$-semigroup property of~$A_\alpha$.

Given a closed operator~$A$ in a Hilbert space,
we denote by $\sigma(A)$, $\sigma_\mathrm{p}(A)$,
$\sigma_\mathrm{d}(A)$, $\sigma_\mathrm{e}(A)$,
repectively, the spectrum, the point spectrum,
the discrete spectrum and the essential spectrum of~$A$.
There are various definitions of the essential spectrum
for non-self-adjoint operators in the literature,
\cf~\cite{Gustafson-Weidmann_1969}.
We define 
$
  \sigma_\mathrm{e}(A)
  := \sigma(A)\setminus\sigma_\mathrm{d}(A)
$ 
and recall that~$\lambda\in\sigma_\mathrm{d}(A)$ 
if and only if it is an isolated point of~$\sigma(A)$
consisting of an eigenvalue with finite algebraic multiplicity
and with the range of $A-\lambda\id$ being closed.  
  
We also point out some special conventions that we use
throughout the paper:
$\Nat^*=\Nat\setminus\{0\}$, where $\Nat=\{0,1,2,\dots\}$,  
and $\Real_\pm^*=(0,\pm\infty)$.

\section{Preliminaries}\label{Sec.Preliminary}
Denote by~$(\cdot,\cdot)$ the inner product in the Hilbert space~$\s^2(\Omega)$
and by~$\|\cdot\|$ the corresponding norm, and let~$\Omega$ be an (arbitrary)
open connected set in $\Real^d$, with $d \geq 1$.

\subsection{The auxiliary Schr\"odinger operator}
%
In this section we consider the family of one-parameter 
Schr\"odinger operators~$\op_\mu$ on~$\s^2(\Omega)$,
subject to Dirichlet boundary conditions,
generated by the differential expression~(\ref{diff.expression}),  
where~$\mu$ is a real parameter
and $a,b: \Omega \to \Real$ are measurable functions
satisfying~$\Hbasic$. 
 
The operators~$\op_\mu$ are introduced as follows.  
Let $\opf_0$ be the sesquilinear form on $\s^2(\Omega)$ defined by
$$
  \opf_0(\phi,\psi) 
  :=  (\nabla\phi,\nabla\psi) + (\phi,b\,\psi) 
  \,,
  \qquad
  \phi,\psi\in\Dom(\opf_0) 
  := \overline{\Smooth_0^\infty(\Omega)}^{\,\|\cdot\|_{\opf_0}}
  \,,
$$
where 
$$
  \|\cdot\|_{\opf_0}^2
  := \|\cdot\|_{\Hilbert^1(\Omega)}^2
  + \| (b-b_\mathrm{min})^\demi \cdot \|^2
  \,.
$$  
(Note that $\Dom(\opf_0)$ is continuously and densely
embedded in $\Hilbert_0^1(\Omega)$.)
Since the form $\opf_0$ is densely defined, closed, symmetric
and bounded from below~\cite[Sec.~VII.1.1]{Edmunds-Evans},
it gives rise to a unique self-adjoint operator~$\op_0$ 
which is also bounded from below. \!%
Using the representation theorem~\cite[Chap.~VI, Thm.~2.1]{Kato},
we have that 
$$
  \op_0\psi = \ope_0\psi \;\!, 
  \qquad 
  \psi\in\Dom(\op_0) = \left\{\psi\in\Dom(\opf_0)\,|\, 
  \ope_0\psi \in \s^2(\Omega) \right\}.
$$ 
Since~$a$ is bounded, the quadratic form
$$
  \opf_\mu[\psi] := \opf_0[\psi]+\mu(\psi,a\psi) \,, 
  \qquad
  \psi\in\Dom(\opf_\mu):=\Dom(\opf_0)
  \,,
$$
gives rise to a self-adjoint operator~$\op_\mu$ 
which is bounded from below and satisfies 
\begin{equation}\label{representation}
  \op_\mu\psi = \ope_\mu\psi \,, 
  \qquad 
  \psi\in\Dom(\op_\mu)=\Dom(\op_0)
  \,.
\end{equation}

Let $\{\gamma_n(\mu)\}_{n\in\Nat^*}$ be the non-decreasing  
sequence of numbers corresponding to the spectral problem of~$\op_\mu$
according to the Rayleigh-Ritz variational formula, \ie\/
\begin{equation}\label{Rayleigh}
  \gamma_n(\mu) := 
  \inf\left\{ \left.
  \sup_{\psi \in L}  
  \frac{\opf_\mu[\psi]}{\,\|\psi\|^2}
  \ \right| \
  L \subseteq \Dom(\opf_0)
  \ \& \ \dim(L)=n
  \right\}.
\end{equation}
\begin{Proposition}\label{Prop.Asymptotes}
Assume~$\Hbasic$.
Each $\mu \mapsto \gamma_n(\mu)$
is a continuous function with the uniform asymptotics
$$
  \gamma_n(\mu) =
\begin{cases}
  a_\mathrm{max} \, \mu + o(\mu)
  & \mbox{as} \quad \mu\to -\infty \,, \\
  a_\mathrm{min} \, \mu + o(\mu)
  & \mbox{as} \quad \mu\to +\infty \,.
\end{cases}
$$ 
\end{Proposition}
\begin{proof}
The continuity is immediate from the definition~(\ref{Rayleigh}); 
indeed, the Lipschitz condition 
$
  |\gamma_n(\mu) - \gamma_n(\mu')| \leq \|a\|_\infty |\mu-\mu'|
$
holds true. 
The asymptotics follow from the more general 
Theorem~\ref{Thm.Asymptotes.App}
in Appendix~\ref{App.Asymptotes}.
\end{proof}
\noindent
We recall that~$\gamma_n(\mu)$ represents either 
a discrete eigenvalue of~$\op_\mu$
or the threshold of its essential spectrum,
\cf~\cite[Sec.~4.5]{Davies};
$\gamma_1(\mu)$~is the spectral threshold of~$\op_\mu$,
\ie, $\gamma_1(\mu)=\inf\sigma(\op_\mu)$.

Let $\gamma_\infty(\mu):=\inf\sigma_\mathrm{e}(\op_\mu)$
denote the threshold of the essential spectrum of~$\op_\mu$
(if $\sigma_\mathrm{e}(\op_\mu)=\varnothing$, 
we put $\gamma_\infty(\mu):=+\infty$). 
The following formula is a generalization (\cf~\cite[Sec.~X.5]{Edmunds-Evans})
of a result due to Persson~\cite{Persson}: 
\begin{equation}\label{Persson}
  \gamma_\infty(\mu) = 
  \sup\left\{ \left.
  \inf_{\psi \in \Smooth_0^\infty(\Omega\setminus K)}  
  \frac{\opf_\mu[\psi]}{\,\|\psi\|^2}
  \ \right| \
  K \ \mbox{compact subset of} \ \Omega 
  \right\}.
\end{equation}
\begin{Proposition}\label{Prop.Asymptotes.ess}
Assume~$\Hbasic$.
$\mu \mapsto \gamma_\infty(\mu)$
is a continuous function satisfying 
$$
  \gamma_\infty(\mu) \geq
\begin{cases}
  \gamma_\infty(0) + a_\mathrm{max} \, \mu 
  & \mbox{if} \quad \mu \leq 0 \,, \\
  \gamma_\infty(0) + a_\mathrm{min} \, \mu 
  & \mbox{if} \quad \mu \geq 0 \,.
\end{cases}
$$ 
\end{Proposition}
\begin{proof}
The Lipschitz condition 
$
  |\gamma_\infty(\mu) - \gamma_\infty(\mu')| \leq \|a\|_\infty |\mu-\mu'|
$
and the lower bounds to~$\gamma_\infty(\mu)$ 
are immediate consequences of~(\ref{Persson}).
\end{proof}
%
 
\subsection{The damped wave operator}
%
Let us introduce the Hilbert space
$
  \Hilbert := \Dom(\opf_0) \times \s^2(\Omega)  
$
of vectors $\Psi\equiv\{\psi_1,\psi_2\}$,
where $\psi_1\in\Dom(\opf_0)$ and $\psi_2\in\s^2(\Omega)$,
and let us equip it with the norm~$\|\cdot\|_\Hilbert$ defined by 
\begin{equation*}
  \|\Psi\|_\Hilbert^2 \equiv \|\{\psi_1,\psi_2\}\|_\Hilbert^2
  := \|\psi_1\|_{\opf_0}^2 + \|\psi_2\|^2   
  \,.
\end{equation*}
It is clear that~$\Hilbert$ is the completion 
of $\Smooth_0^\infty(\Omega)\times\Smooth_0^\infty(\Omega)$
w.r.t.~$\|\cdot\|_\Hilbert$.

Writing $\psi_1:=\psi$ and $\psi_2:=\psi_{t}$, 
the problem~(\ref{wave.eq}) is equivalent to the first-order system
\begin{equation}\label{wave.eq.bis}
\left\{ 
\begin{aligned}
  \Psi_{t} &= A_\alpha \Psi  
  \,, 
  \\
  \Psi(0,\cdot) & = \Phi \in \Dom(A_\alpha)
  \,,   
\end{aligned}
\right.
\end{equation}
where~$A_\alpha$ is the operator associated with~(\ref{wave.op.formal});
it is properly defined as follows:
\begin{align}\label{wave.op}
  A_\alpha \{\psi_1,\psi_2\} 
  & := \big\{\psi_2,-\op_0\psi_1-\alpha a\psi_2\big\}
  \,,
  \\
  \{\psi_1,\psi_2\} \in \Dom(A_\alpha) 
  & := \Dom(\op_0)\times\Dom(\opf_0)
  \,. \nonumber
\end{align}
\begin{Proposition}\label{Prop.closed}
Assume~$\Hbasic$.
The operator~$A_\alpha$ is densely defined and closed.
Furthermore, it is the generator of a $C_0$-semigroup on~$\Hilbert$.
\end{Proposition}
\begin{proof}
The density and closedness follow due to 
the same properties for~$\op_0$, $\opf_0$
and the boundedness assumption on~$a$.  
The $C_0$-semigroup property is proved in Appendix~\ref{App.C0}
(\cf~Theorem~\ref{Thm.C0}). 
\end{proof}

We shall use the following sufficient condition
(characterization of the essential spectrum due to Wolf~\cite{Wolf}).
\begin{Lemma}\label{Lem.Weyl}
Assume~$\Hbasic$.
If there exists a singular sequence of~$A_\alpha$
corresponding to $\lambda\in\Com$, \ie,
a sequence $\{\Psi_n\}_{n\in\Nat}\subseteq\Dom(A_\alpha)$
with the properties
$$
  \|\Psi_n\|_\Hilbert = 1 \,, 
  \qquad
  \Psi_n \xrightarrow[n\to\infty]{w} 0 
  \quad\mbox{and}\quad
  (A_\alpha-\lambda\id)\Psi_n \xrightarrow[n\to\infty]{} 0 
  \quad\mbox{in}\quad
  \Hilbert
  \,,
$$
then $\lambda\in\sigma_\mathrm{e}(A_\alpha)$.
\end{Lemma}

We recall that the existence of singular sequence  
provides also a necessary condition in the characterization
of essential spectrum for self-adjoint operators.
 
\section{Real spectrum of the damped wave operator}\label{Sec.Real} 
%
Our further approach is based on the following trivial,
but essential, observation that links the spectral problems 
for the operators~$A_\alpha$ and~$\op_\mu$. 
\begin{Lemma}\label{Lem.basic}
Assume~$\Hbasic$.
$\forall\mu\in\Real$, \ $\forall\alpha>0$,  
\begin{align*}
  \emph{(i)} && 
  - (\mu/\alpha)^2 &\in \sigma_\mathrm{p}(\op_\mu) 
  & \rule{-5ex}{0ex} \Longleftrightarrow\qquad
  \mu/\alpha &\in \sigma_\mathrm{p}(A_\alpha) \,,
  \\
  \emph{(ii)} && 
  - (\mu/\alpha)^2 &\in \sigma_\mathrm{e}(\op_\mu) 
  & \rule{-5ex}{0ex} \Longrightarrow\qquad
  \mu/\alpha &\in \sigma_\mathrm{e}(A_\alpha)  \,.
\end{align*}
\end{Lemma}
\begin{proof}
\ad~(i). 
If $- (\mu/\alpha)^2 \in \sigma_\mathrm{p}(\op_\mu)$,
then there exists a $\psi\in\Dom(\op_0)$ such that 
$\op_\mu\psi=- (\mu/\alpha)^2\psi$ and 
$
  \Psi\equiv\{\psi_1,\psi_2\}:=\{\psi,(\mu/\alpha)\psi\}
  \in \Dom(A_\alpha)
$
is easily checked to satisfy
$
  A_\alpha \Psi = (\mu/\alpha) \Psi
$. 
Conversely, if~$\mu/\alpha\in\sigma_\mathrm{p}(A_\alpha)$,
then there exist $\psi_1\in\Dom(\op_0)$ and $\psi_2\in\Dom(\opf_0)$
satisfying the system
$
  \psi_2 = (\mu/\alpha) \psi_1 \,,
  -\op_0\psi_1-\alpha a\psi_2 = (\mu/\alpha) \psi_2 \,,
$
which yields $\op_\mu\psi_1=- (\mu/\alpha)^2\psi_1$. 

\ad~(ii).
If $- (\mu/\alpha)^2 \in \sigma_\mathrm{e}(\op_\mu)$, 
then there exists a sequence 
$\{\psi_n\}_{n\in\Nat}\subseteq\Dom(\op_0)$
for which $\|\psi_n\| = 1$, 
$\psi_n \xrightarrow[]{w} 0$
and $(\op_\mu+(\mu/\alpha)^2\id)\psi_n \xrightarrow[]{} 0$
in $\s^2(\Omega)$ as~$n\to\infty$.
We shall show that the sequence 
$
  \{\Psi_n/\|\Psi_n\|_\Hilbert\}_{n\in\Nat}
$,
where 
$
  \Psi_n 
  \equiv \{\psi_n^1,\psi_n^2\}
  :=\{\psi_n,(\mu/\alpha)\psi_n\}
$,
is the singular sequence of Lemma~\ref{Lem.Weyl}
with~$\lambda=\mu/\alpha$.
First of all, we note that
$
  \|\Psi_n\|_\Hilbert
  \geq \|\psi_n\|_{\opf_0}
  \geq \|\psi_n\|
  = 1
$,
so the sequence is well-defined, normalized in~$\Hilbert$,
and it is enough to check the weak and strong convergence
of Lemma~\ref{Lem.Weyl} for the sequence $\{\Psi_n\}_{n\in\Nat}$.
Clearly, 
$$
  (A_\alpha - (\mu/\alpha)\id) \Psi_n
  = \{0,-(\op_\mu+(\mu/\alpha)^2\id)\psi_n\}
  \xrightarrow[n\to\infty]{} \{0,0\}
  \quad\mbox{in}\quad
  \Dom(\opf_0)\times\s^2(\Omega)
  .
$$ 
It is also clear that $\psi_n^2 \xrightarrow[]{w} 0$ 
in~$\s^2(\Omega)$ as~$n\to\infty$,
so it remains to show that $\psi_n^1 \equiv \psi_n \xrightarrow[]{w} 0$ 
in~$\Dom(\opf_0)$ as~$n\to\infty$.
The latter can be seen by writing
$$
  (\phi,\psi_n)_{\opf_0}
  = \left( \phi, [\op_\mu+(\mu/\alpha)^2\id] \psi_n \right)
  + \left(\phi,
  [1 - b_\mathrm{min} - \mu a - (\mu/\alpha)^2] \psi_n 
  \right) 
$$
for every $\phi\in\Dom(\opf_0)$.
\end{proof}

The importance of this result is that it makes it possible
to look for real points in the spectrum of~$A_\alpha$
by considering the much simpler (since self-adjoint) 
spectral problem for~$\op_\mu$.
(Although the latter continues to make sense for complex~$\mu$,
in such a case, however, one is merely trading
one non-self-adjoint problem for another.) 
In view of our definition of essential spectrum,
it is clear that, in general, there might be a 
$\mu/\alpha \in \sigma_\mathrm{e}(A_\alpha)$
such that $- (\mu/\alpha)^2 \not\in \sigma_\mathrm{e}(\op_\mu)$,  
and that is why we do not have a proof 
for the converse implication in~(ii) of Lemma~\ref{Lem.basic}.  
  
Lemma~\ref{Lem.basic} yields 
\begin{Proposition}\label{Prop.basic}
Assume~$\Hbasic$.
For any~$\alpha>0$,
\begin{align*}
  \emph{(i)} && 
  0\in\sigma_\mathrm{p}(\op_0)  
  \quad\Longleftrightarrow\quad &  
  0\in\sigma_\mathrm{p}(A_\alpha) \,,
  \\
  \emph{(ii)} && 
  0\in\sigma_\mathrm{e}(\op_0)  
  \quad\Longrightarrow\quad &  
  0\in\sigma_\mathrm{e}(A_\alpha) \,,
  \\ 
  \emph{(iii)} && 
  \gamma_1(0) < 0
  \quad\Longrightarrow\quad &  
  \left(
  \sigma(A_\alpha) \cap \Real_-^* \not= \varnothing
  \quad\&\quad
  \sigma(A_\alpha) \cap \Real_+^* \not= \varnothing
  \right) \,,
  \\ 
  \emph{(iv)} && 
  \gamma_\infty(0) < 0 
  \quad\Longrightarrow\quad &  
  \left(
  \sigma_\mathrm{e}(A_\alpha) \cap \Real_-^* \not= \varnothing
  \quad\&\quad
  \sigma_\mathrm{e}(A_\alpha) \cap \Real_+^* \not= \varnothing
  \right) \,.
\end{align*}
\end{Proposition}
\begin{proof}
(i) and~(ii) are direct consequences of Lemma~\ref{Lem.basic}.
Since $\gamma_1(0)<0$, it follows by Proposition~\ref{Prop.Asymptotes}
that there exist at least one negative and one positive~$\mu$
at which the curve $\mu\mapsto\gamma_1(\mu)\in\sigma(\op_\mu)$ 
intersects the parabola $\mu\mapsto -(\mu/\alpha)^2$.
If $\gamma_\infty(0)<0$, one applies the same argument to~$\gamma_\infty$,
with help of Proposition~\ref{Prop.Asymptotes.ess}.
\end{proof}

In~(iii) of Proposition~\ref{Prop.basic},
since the spectral nature of~$\gamma_1(\mu)$,
\ie, whether it belongs to the discrete or essential spectrum of~$\op_\mu$,
may vary with~$\mu$, we do not know of which spectral nature
are the corresponding points in the real spectrum of~$A_\alpha$. 
This can be decided, however, if~$\op_\mu$ has no negative essential spectrum.
Let us define
$$
  N_\mu := \#\big\{
  n\in\Nat^* \,\big| \ 
  \gamma_n(\mu) < \min\{0,\gamma_\infty(\mu)\}
  \big\}
  \,,
$$
\ie, the number of negative eigenvalues of~$\op_\mu$
below its essential spectrum, counting multiplicities.
\begin{Proposition}\label{Prop.basic.evs}
Assume~$\Hbasic$.
If $\gamma_\infty(\mu) \geq 0$ for all $\mu \leq 0$
(respectively all \mbox{$\mu \geq 0$}),
then there exist at least~$N_0$ negative
(respectively~$N_0$ positive) eigenvalues of~$A_\alpha$
for all~$\alpha>0$. 
\end{Proposition}
\begin{proof}
The claim is trivial for~$N_0=0$. 
Let us assume $N_0 \geq 1$ and $\gamma_\infty(\mu) \geq 0$ 
for all $\mu \leq 0$. Under the assumptions,  
$\gamma_n(0)\in\sigma_\mathrm{d}(\op_0)\cap\Real_-^*$, $n\in\{1,\dots,N_0\}$,
and one applies the argument of the proof of Proposition~\ref{Prop.basic}
to every curve $\mu\mapsto\gamma_n(\mu)$.
The difference is now that~$\gamma_n(\mu)$
cannot become~$\gamma_\infty(\mu)$ for $\mu \leq 0$, 
without crossing the parabola $\mu\mapsto -(\mu/\alpha)^2$,
so the intersection corresponds to an eigenvalue.
The same argument holds for the case~$\mu \geq 0$.  
\end{proof}

Proposition~\ref{Prop.basic.evs} is a generalization 
of~\cite[Prop.~3.2]{Freitas_1996}, where the result 
was established for both bounded~$\Omega$ and~$b$;
then $\gamma_\infty(\mu) = +\infty$ 
(\ie, $\sigma_\mathrm{e}(\op_\mu)=\varnothing$)
and $N_0<+\infty$.  
It is worth to notice that, in our more general setting,
$\gamma_\infty(\mu)$ can be finite and~$N_0$ infinite. 
 
To state and prove the main result of the paper,
we impose on the damping~$a$ the following indefinitness condition
\begin{Assumption}{\Hsign}
$a_\mathrm{min}<0$ \ and \ $a_\mathrm{max}>0$. 
\end{Assumption}
\begin{Theorem}\label{Thm.main}
Assume~$\Hbasic$ and~$\Hsign$.
For all sufficiently large~$\alpha$,
there is at least one negative and one positive point
in the spectrum of~$A_\alpha$. 
Moreover,
\begin{itemize}
\item[\emph{(i)}]
if 
$
  {\displaystyle
  \lim_{\mu\to-\infty}\gamma_\infty(\mu)/\mu 
  < a_\mathrm{max}
  }
$
(respectively  
$
  {\displaystyle
  \lim_{\mu\to+\infty}\gamma_\infty(\mu)/\mu 
  > a_\mathrm{min}
  }
$),
then there exists a positive increasing sequence
$\{\alpha_j^-\}_{j=1}^\infty$
(respectively $\{\alpha_j^+\}_{j=1}^\infty$) such that: 
$\forall\alpha\in[\alpha_j^-,\alpha_{j+1}^-)$
there exist at least $j$ negative eigenvalues  of~$A_\alpha$
(respectively
$\forall\alpha\in[\alpha_j^+,\alpha_{j+1}^+)$
there exist at least $j$ positive eigenvalues  of~$A_\alpha$);
\item[\emph{(ii)}]
if $\gamma_\infty(\mu) > 0$ for all $\mu \leq 0$
(respectively for all $\mu \geq 0$),
then there exists a positive increasing sequence
$\{\alpha_j^-\}_{j=1}^\infty$
(respectively $\{\alpha_j^+\}_{j=1}^\infty$) such that:
$\forall\alpha\in(\alpha_j^-,\alpha_{j+1}^-]$
there exist at least $2 j + N_0$ negative eigenvalues  of~$A_\alpha$
(respectively $\forall\alpha\in(\alpha_j^+,\alpha_{j+1}^+]$
there exist at least $2 j + N_0$ positive eigenvalues  of~$A_\alpha$)
and $\forall\alpha\in(0,\alpha_1^-]$
there exist at least~$N_0$ negative eigenvalues of~$A_\alpha$
(respectively $\forall\alpha\in(0,\alpha_1^+]$
there exist at least~$N_0$ positive eigenvalues of~$A_\alpha$);
\item[\emph{(iii)}]
if $\gamma_\infty(\mu) < 0$ 
for some~$\mu \leq 0$ (respectively for some~$\mu \geq 0$),
then there is at least one negative (respectively one positive) 
point in~$\sigma_\mathrm{e}(A_\alpha)$ 
for all sufficiently large~$\alpha$;
\item[\emph{(iv)}]
if $\gamma_\infty(0) > 0$ and $\gamma_\infty(\mu) < 0$ 
for some~$\mu \leq 0$ (respectively for some $\mu \geq 0$),
then there are at least two negative (respectively two positive) 
points in~$\sigma_\mathrm{e}(A_\alpha)$ 
for all sufficiently large~$\alpha$.
\end{itemize}
\end{Theorem}
\begin{proof}
We shall prove the claims for the case~$\mu \geq 0$
(\ie, the results about the positive spectrum of~$A_\alpha$),
the case of~$\mu \leq 0$ being similar.
 
If~$\gamma_1(0) < 0$, then the first statement holds 
for all~$\alpha>0$ by Proposition~\ref{Prop.basic}.
If~$\gamma_1(0) > 0$, then the curve $\mu\mapsto\gamma_1(\mu)$
does not intersect the parabola $\mu\mapsto -(\mu/\alpha)^2$
for sufficiently small~$\alpha$, 
but it does (in fact, at least twice),
in view of Proposition~\ref{Prop.Asymptotes} and~$\Hsign$,
if~$\alpha$ is large enough;
the intersections are positive and the result follows
by Lemma~\ref{Lem.basic}.
If~$\gamma_1(0) = 0$, then there will be at least 
one positive intersection by the same argument.
 
\ad~(i).
By virtue of Proposition~\ref{Prop.Asymptotes}, the assumption~$\Hsign$
and the condition about the asymptotic behaviour 
of the threshold of the essential spectrum,
it follows that for any $n\in\Nat^*$ there is an $\mu_n>0$ such that 
$\gamma_n(\mu) < \min\{0,\gamma_\infty(\mu)\}$ for all $\mu\geq\mu_n$,
\ie, $\gamma_n(\mu)$ is a negative discrete eigenvalue of~$\op_\mu$
for all $\mu\geq\mu_n$. 
Moreover, the sequence~$\{\mu_n\}_{n=1}^\infty$
can be chosen in such a way that 
the sequence~$\{\alpha_{n}^+\}_{n=1}^\infty$ 
defined by $\alpha_{n}^+:=\mu_n/\sqrt{-\gamma_n(\mu_n)}$
is increasing.
However, by Proposition~\ref{Prop.Asymptotes},
each curve $\mu\mapsto\gamma_n(\mu)$ has at least 
one intersection with the parabola $\mu\mapsto -(\mu/\alpha)^2$
for all $\alpha\in[\alpha_{n}^+,\infty)$. 
Hence, by Lemma~\ref{Lem.basic}, 
$A_\alpha$ has at least~$n$ positive eigenvalues
for all 
$
  \alpha\in[\alpha_{n}^+,\alpha_{n+1}^+)
$.
 
\ad~(ii).
The~$N_0$ positive eigenvalues of~$A_\alpha$ correspond,
by Proposition~\ref{Prop.basic.evs},
to~$N_0$ negative eigenvalues of~$\op_0$, 
for all~$\alpha>0$.  
Let~$n_0>N_0$ be the lowest of all~$n$ satisfying $\gamma_n(0)>0$
(it always exists since $\gamma_\infty(0) > 0$). 
Let~$n\in\{n_0,\dots,\infty\}$. 
Since $\gamma_n(0)>0$, the curve $\mu\mapsto\gamma_n(\mu)$
does not intersect the parabola $\mu\mapsto -(\mu/\alpha)^2$
for~$\alpha$ sufficiently small. 
It follows by Proposition~\ref{Prop.Asymptotes} that,
by increasing~$\alpha$, it is possible to make 
these two curves touch for some~$\mu>0$. 
Denote by~$\tilde{\alpha}_{n-n_0+1}^+$ the point 
for which this happens for the first time.
Then, as~$\alpha$ increases past this value,
there will always exist at least two intersections for positive~$\mu$.
Hence, by Lemma~\ref{Lem.basic}, 
$A_\alpha$ has at least $2(n-n_0+1)$ positive eigenvalues
for all 
$
  \alpha \in (\tilde{\alpha}_{n-n_0+1}^+,\infty)
$.
Since~$\gamma_n$ is below~$\gamma_{n+1}$,
it follows that the sequence~$\{\tilde{\alpha}_j^+\}_{j=1}^\infty$ 
is not decreasing.
In fact, there exists an increasing subsequence
$\{\alpha_j^+\}_{j=1}^\infty$ because each negative~$\gamma_n(\mu)$
represents an eigenvalue of \emph{finite} multiplicity.

\ad~(iii).
Given a~$\mu \geq 0$ such that $\gamma_\infty(\mu)<0$,
we define $\alpha_0:=\mu/\sqrt{-\gamma_\infty(\mu)}$
and it follows by Proposition~\ref{Prop.Asymptotes.ess}
that there is at least one positive intersection 
of $\mu\mapsto\gamma_\infty(\mu)$ 
with the parabola $\mu\mapsto -(\mu/\alpha)^2$
for every $\alpha \in [\alpha_0,\infty)$.   
The result then follows by Lemma~\ref{Lem.basic}.

\ad~(iv).
Since $\gamma_\infty(0)>0$, the curve $\mu\mapsto\gamma_\infty(\mu)$
does not intersect the parabola $\mu\mapsto -(\mu/\alpha)^2$
for~$\alpha$ sufficiently small. 
However, since there is a $\mu>0$ such that $\gamma_\infty(\mu)<0$,
it follows that, by increasing~$\alpha$, it is possible to make 
these two curves touch for some~$\mu>0$. 
Denote by~$\alpha_0$ the point 
for which this happens for the first time.
Then, by Proposition~\ref{Prop.Asymptotes.ess},
as~$\alpha$ increases past this value,
there will always exist at least two intersections for positive~$\mu$.
Hence, by Lemma~\ref{Lem.basic}, 
$A_\alpha$ has at least two positive points
in the its essential spectrum
for all 
$
  \alpha \in (\alpha_0,\infty)
$.
\end{proof} 
\begin{Remarks}
Assume~$\Hsign$ and $\gamma_\infty(\mu)>0$ for all~$\mu\in\Real$,
and denote by~$\lambda(\alpha)$ a real eigenvalue of~$A_\alpha$,
as described in the proof of the part~(ii) of Theorem~\ref{Thm.main}.
Following~\cite{Freitas_1996}, 
it is possible to give some more properties of~$\lambda(\alpha)$. 
 
For instance, the following properties are clear 
from the proof of Proposition~\ref{Prop.basic.evs}
and the part~(ii) of Theorem~\ref{Thm.main}
(\cf~also~\cite[proof of Prop.~3.8]{Freitas_1996}):
If $\lambda(\alpha)<0$
(respectively $\lambda(\alpha)>0$) corresponds to a negative
eigenvalue of~$\op_0$, the function $\alpha\mapsto\lambda(\alpha)$
is decreasing (respectively increasing) and 
$\lambda(\alpha) \to -\infty$ 
(respectively~$+\infty$) as $\alpha\to+\infty$.
In the case of the two negative (respectively positive) eigenvalues,
corresponding to some positive~$\gamma_n(0)$,
one of these eigenvalues is increasing (respectively decreasing) to zero,
while the other is decreasing to~$-\infty$ 
(respectively increasing to~$+\infty$). 

Also, since~$\gamma_1$ remains below all other curves~$\gamma_n$ 
and the eigenfunction of~$\op_\mu$ corresponding to~$\gamma_1(\mu)<0$
can be chosen to be positive, we get the following characterization
of the eigenfunctions of~$A_\alpha$ corresponding 
to its extreme real eigenvalues 
(\cf~\cite[Thm.~3.9]{Freitas_1996}):
If~$N_0>0$, then the eigenfunction~$\Psi\equiv\{\psi_1,\psi_2\}$ 
corresponding to the smallest negative (respectively largest positive) 
eigenvalue of~$A_\alpha$ can be chosen in such a way that
$\psi_1(x)>0$ and $\psi_2(x)<0$ 
(respectively $\psi_2(x)>0$) for all~$x\in\Omega$.
If~$N_0=0$ and~$A_\alpha$ has negative (respectively positive) eigenvalues
(this always happens for~$\alpha$ sufficiently large),
then the eigenfunctions~$\Psi\equiv\{\psi_1,\psi_2\}$ 
corresponding to the largest or the smallest of these eigenvalues
can be chosen in such a way that $\psi_1(x)>0$ and $\psi_2(x)<0$ 
(respectively $\psi_2(x)>0$) for all~$x\in\Omega$.
\end{Remarks}

If both~$\Omega$ and~$b$ are bounded,
Theorem~\ref{Thm.main} reduces to~\cite[Thm.~3.6]{Freitas_1996}.
Actually, all the complexities of the present theorem 
are due to the possible presence of essential spectrum
in our more general setting. 
The question of essential spectrum is discussed
in more details in the following section.
 
\section{More on the essential spectrum}\label{Sec.Ess}
%
All the spectral results of the previous section
are based on proving an intersection of the functions~$\gamma_n$
or~$\gamma_\infty$ with a parabola. 
The results are qualitative since we have used just the continuity
and asymptotic properties of the functions
(Propositions~\ref{Prop.Asymptotes} and~\ref{Prop.Asymptotes.ess}).
The purpose of this section is twofold.
Firstly, we use known properties of~$\gamma_\infty$
in the case where~$a$ converges to a fixed value at infinity
in order to state more precise results 
about the real essential spectrum of~$A_\alpha$.
Secondly, on characteristic examples,
we discuss the nature of the real spectrum of~$A_\alpha$ 
as related to the form of the domain~$\Omega$
or the behaviour of the potential~$b$.
 
\subsection{Asymptotically constant damping} 
%
If~$\gamma_n(\mu)$ is a discrete eigenvalue of~$\op_\mu$
for all~$\mu\in\Real$ and~$a$ is not constant, 
then it is hard to say anything about
the behaviour of the curve $\mu\mapsto\gamma_n(\mu)$,
apart from small values of the parameter~$\mu$ 
(via perturbation techniques) 
or large~$\mu$ (\cf~Proposition~\ref{Prop.Asymptotes}).
On the other hand, the behaviour of the curve 
$\mu\mapsto\gamma_\infty(\mu)$ is expected to be usually simple
because, heuristically speaking, the essential spectrum comes
from ``what happens very far away" only.
The above statement is made precise
in context of the following result, which can be established 
by means of the decomposition principle~\cite[Chap.~X]{Edmunds-Evans}. 
\begin{Proposition}\label{Prop.decomposition}
Assume~$\Hbasic$. 
If~$\Omega$ is unbounded and the limit
$$
  a_\infty := \lim_{|x|\to\infty,\,x\in\Omega} a(x) 
$$
exists, then
$$
  \forall \mu\in\Real,
  \qquad
  \sigma_\mathrm{e}(\op_\mu) 
  = \sigma_\mathrm{e}(\op_0) + a_\infty \, \mu
  \,.
$$
\end{Proposition}

In particular, under the assumptions of Proposition~\ref{Prop.decomposition},
the curve $\mu\mapsto\gamma_\infty(\mu)$ is linear 
(if $\gamma_\infty(0)=+\infty$,
then $\gamma_\infty(\mu)=+\infty$ for all~$\mu\in\Real$).
Applications of this result to Theorem~\ref{Thm.main} are obvious.
For instance, the sufficient conditions 
from the case~(i) of Theorem~\ref{Thm.main}
are satisfied if $a_\infty<a_\mathrm{max}$
(respectively $a_\infty>a_\mathrm{min}$).
Also, it follows from the case~(iii) of Theorem~\ref{Thm.main}
that there is always a negative (respectively positive), 
point in $\sigma_\mathrm{e}(A_\alpha)$
provided $\sigma_\mathrm{e}(\op_0)\not=\varnothing$
and $a_\infty>0$ (respectively $a_\infty<0$).
Actually, a stronger result holds if the essential spectrum of~$\op_0$  
is an infinite interval:
\begin{Proposition}\label{Prop.Ess}
Assume the hypotheses of \emph{Proposition~\ref{Prop.decomposition}}.
Suppose also that
$
  \sigma_\mathrm{e}(\op_0) = [\gamma_\infty(0),\infty)
$
and
$
  \delta_\alpha
  := \alpha^2 a_\infty^2 - 4 \gamma_\infty(0) 
  \geq 0
$.
Then
$$
  \left[
  \mbox{$\demi$}\left(-\alpha a_\infty-\sqrt{\delta_\alpha}\right),
  \mbox{$\demi$}\left(-\alpha a_\infty+\sqrt{\delta_\alpha}\right)
  \right]
  \subseteq \sigma_\mathrm{e}(A_\alpha)
  \,.
$$
\end{Proposition}
\begin{proof}
Combining the assumption 
$
  \sigma_\mathrm{e}(\op_0) = [\gamma_\infty(0),\infty)
$,
Proposition~\ref{Prop.decomposition} and Lemma~\ref{Lem.basic},
we see that every point between the intersections (if they exist) 
of the parabola $\mu\mapsto -(\mu/\alpha)^2$
with the line $\mu\mapsto\gamma_\infty(0)+a_\infty\mu$
lies in the essential spectrum of~$A_\alpha$.
The claim then follows by solving the corresponding quadratic equation
(in particular, the condition $\delta_\alpha \geq 0$
ensures its solvability).
\end{proof} 
%

\subsection{Examples} 
%
In this subsection the significant features of~$\Omega$ or~$b$
as regards $\sigma_\mathrm{e}(A_\alpha)$ are discussed.
Since the spectrum of~$A_\alpha$ is purely discrete 
whenever~$\Omega$ is bounded, we restrict ourselves 
to unbounded~$\Omega$ here.
We also assume that the hypotheses of Proposition~\ref{Prop.decomposition}
hold true, in particular, $a_\infty$~denotes the asymptotic 
value of the damping~$a$. 

\subsubsection{Different types of domains} 
%
We adopt the following classification of domains  
in Euclidean space due to Glazman~\cite[\S\,49]{Glazman} 
(see also~\cite[Sec.~X.6.1]{Edmunds-Evans}):
$\Omega$ is \emph{quasi-conical} 
if it contains arbitrarily large balls;
$\Omega$ is \emph{quasi-cylindrical}
if it is not quasi-conical but contains a sequence 
of identical disjoint balls;
$\Omega$ is \emph{quasi-bounded}
if it is neither quasi-conical nor quasi-cylindrical.
For simplicity, let us assume that~$b=0$,
so that~$\op_0$ is just the Dirichlet Laplacian in~$\Omega$.  

\paragraph{Quasi-bounded domains} 
%
If the boundary of~$\Omega$ is sufficiently smooth
(\cf~\cite[Sec.~X.6.1]{Edmunds-Evans} for precise 
conditions in terms of capacity), 
Mol\-\v{c}anov's criterion~\cite{Molcanov} implies that 
the spectrum of~$\op_0$ is purely discrete.
Then the same is true for the operator~$A_\alpha$, \ie,
$
  \sigma_\mathrm{e}(A_\alpha) = \varnothing
$.
If the hypothesis~$\Hsign$ holds true,
the part~(ii) of Theorem~\ref{Thm.main} claims
that there are always negative and positive 
eigenvalues of~$A_\alpha$ for~$\alpha$ sufficiently large.
This situation is analogous to the bounded case~\cite{Freitas_1996}. 

\paragraph{Quasi-conical domains} 
%
The spectrum of~$\op_0$ is purely essential 
and equals the interval~$[0,\infty)$. 
By Propositions~\ref{Prop.decomposition} and~\ref{Prop.Asymptotes},
$\sigma_\mathrm{e}(\op_\mu)=[a_\infty\mu,\infty)$,
and there are discrete eigenvalues below the essential spectrum 
for sufficiently large positive~$\mu$ if $a_\infty>a_\mathrm{min}$  
or sufficiently large negative~$\mu$ if $a_\infty<a_\mathrm{max}$. 
By Proposition~\ref{Prop.Ess}, 
$$
  \left[
  -\mbox{$\frac{\alpha}{2}$}
  \left(|a_\infty|+a_\infty\right),
  \mbox{$\frac{\alpha}{2}$}
  \left(|a_\infty|-a_\infty\right)
  \right]
  \subseteq \sigma_\mathrm{e}(A_\alpha)
  \,,
$$
\ie, $A_\alpha$ has a real essential spectrum
for any $\alpha \geq 0$.
Information about the real point spectrum of~$A_\alpha$
is contained in the part~(i) of Theorem~\ref{Thm.main}.

\paragraph{Quasi-cylindrical domains} 
%
Since the precise location of $\sigma_\mathrm{e}(\op_0)$ 
is difficult in this situation (a detailed analysis of this problem
may be found in Glazman's book~\cite{Glazman}),
we shall rather illustrate the nature of~$\sigma(A_\alpha)$
in the particular case of tubes.

Let~$\Omega$ be a tubular neighbourhood of radius~$r>0$
about an infinite curve in~$\Real^d$, with $d \geq 2$,
and assume that it does not overlap itself 
and that the condition $\|\kappa_1\|_\infty \, r < 1$ holds true,
where~$\kappa_1$ denotes the first curvature of the reference curve. 
If the tube is asymptotically straight in the sense
that~$\kappa_1$ vanishes at infinity, it was shown in~\cite{ChDFK} 
that the essential spectrum of~$\op_0$ is the interval $[\nu_1,\infty)$,
where~$\nu_1>0$ denotes the first eigenvalue of the Dirichlet Laplacian  
in the $(d-1)$-dimensional ball of radius~$r$. 
(Furthermore, there exists at least one discrete eigenvalue
in~$(0,\nu_1)$ whenever~$\kappa_1\not=0$.) 
It follows that the assumptions of Proposition~\ref{Prop.Ess}
are satisfied if~$a_\infty\not=0$ 
and $\alpha$ is sufficiently large. 
That is, depending on the sign of~$a_\infty$,
$A_\alpha$ always has an open interval 
in its negative or positive essential spectrum 
for~$\alpha$ sufficiently large. 
An information about the real point spectrum of~$A_\alpha$
is contained in the part~(i) or~(ii) of Theorem~\ref{Thm.main}.

\subsubsection{The behaviour of the potential} 
%
For simplicity, let us assume that~$\Omega$ 
is the whole space~$\Real^d$ (a special case of quasi-conical domains). 
Since the potential~$b$ is assumed to be locally bounded,
the essential spectrum of~$\op_0$ 
will depend on the behaviour of~$b$ at infinity only. 
Assuming
$$
  \lim_{|x|\to\infty} b(x) = b_\infty \in [b_\mathrm{min},+\infty]
  \,,
$$
we distinguish two particular behaviours: 

\paragraph{\underline{$b_\infty=+\infty$}} 
%
It is easy to see that $\gamma_\infty(0)=+\infty$
and the spectrum of~$\op_0$ is therefore purely discrete. 
Consequently, the spectrum of~$A_\alpha$ is purely discrete as well
and the situation is analogous with the bounded case~\cite{Freitas_1996}. 

\paragraph{\underline{$b_\infty<+\infty$}} 
%
The essential spectrum of~$\op_0$ is the interval $[b_\infty,\infty)$
and there might be discrete eigenvalues in $[b_\mathrm{min},b_\infty)$
if $b_\mathrm{min}<b_\infty$.
By Proposition~\ref{Prop.decomposition}, we have
$\sigma_\mathrm{e}(\op_\mu)=[b_\infty+a_\infty\mu,\infty)$,
and this gives an information about 
the real essential (respectively real point) 
spectrum of~$A_\alpha$ through Proposition~\ref{Prop.Ess}
(respectively through the part~(i) or~(ii) of Theorem~\ref{Thm.main}).

\section{Remarks and open questions\label{remarks}}
%
In a similar way to what was done in~\cite{Freitas_1996}, the results given
here can be used to obtain information about semilinear wave equations of the
form 
\begin{equation}\label{semilin}
u_{tt} + f(x,u,u_{t}) = \Delta u, \;\; x\in\Omega.
\end{equation}
In particular, the stability of stationary solutions 
of~(\ref{semilin}) can be related to that of solutions of the parabolic equation
\begin{equation}\label{parab}
    u_{t}+f(x,u,0) = \Delta u, \;\; x\in\Omega.
\end{equation}
A typical example of such a result would be the following, the proof 
of which follows in the same way as that of the corresponding result 
in~\cite{Freitas_1996} -- note that equations~(\ref{semilin}) 
and~(\ref{parab}) have the same set of stationary solutions.
\begin{Theorem}
    Let $w$ be a stationary solution of equation~\emph{(\ref{semilin})}. Then 
    if $w$ is linearly unstable when considered as a stationary solution of 
    the parabolic equation~\emph{(\ref{parab})}, it is also linearly unstable 
    as a stationary solution of~\emph{(\ref{semilin})}.
\end{Theorem}
More specific results in the spirit of those in~\cite{Freitas_1996} 
may be given in the case of quasi-bounded domains, for instance.

Regarding stability issues, and as in the bounded case, it would be of interest
to know if there are 
situations when the trivial solution of~(\ref{wave.eq}) remains 
(asymptotically) stable for small $\alpha$ when $a$ is allowed to become
negative. In the former case, the problem remains open in dimension 
higher than one, but there are several results for the one 
dimensional problem showing that this is still 
possible~\cite{bera,FZ}. However, even in the case of the real line, the
situation for unbounded domains seems to be more difficult due to the presence
of essential 
spectrum. Also, the proof of stability in the case of an interval 
used information about the asymptotic behaviour of the spectrum which is 
not available in the unbounded case.

Finally, let us remark that it is possible to adapt 
the approach of the present paper to damped wave equations
with a more general elliptic operator 
instead of~$S_0$ considered here.

\appendix
%
\section{Eigenvalue asymptotics}\label{App.Asymptotes}
%
In this appendix we prove the asymptotics 
of Proposition~\ref{Prop.Asymptotes}.
We proceed in a greater generality by admitting
local singularities of the function~$b$
and that~$a$ is bounded just from below. 
That is, the assumption~$\Hbasic$ is replaced by
\begin{Assumption}{\Hbasicbis}   
$a,b\in\sloc^1(\Omega)$ and
\begin{itemize}
\item[\emph{(a)}]
$a_\mathrm{min} := \essinf a > -\infty$,
\item[\emph{(b)}]
$b=b_1+b_2$ where 
$b_1\in\sloc^1(\Omega)$,
$b_{1,\mathrm{min}} := \essinf b_1 > -\infty$ \\
and $b_2\in\s^p(\Omega)$ for some $p \in [d,\infty]$. 
\end{itemize}
\end{Assumption}

Here we restrict ourselves to~$\mu \geq 0$
(the case of~$\mu \leq 0$ being included if~$a\in\s^\infty(\Omega)$)
and define the realization~$\op_\mu$  
of the differential expression~(\ref{diff.expression})
as the operator associated with the sesquilinear form 
\begin{align*}
  \opf_\mu(\phi,\psi) 
  &:= 
  (\nabla\phi,\nabla\psi) + \mu \, (\phi,a\psi) + (\phi,b\psi) 
  \,,
  \\
  \phi,\psi\in\Dom(\opf_\mu) 
  &:= \overline{\Smooth_0^\infty(\Omega)}^{\,\|\cdot\|_{\opf_\mu}}
  \,,
  \\
  \|\cdot\|_{\opf_\mu}^2
  &:= \|\cdot\|_{\Hilbert^1(\Omega)}^2
  + \big\|
  \big[\mu \, (a-a_\mathrm{min})+b_1-b_{1,\mathrm{min}}+b_2^+\big]^\demi
  \cdot\big\|^2
  \,,
\end{align*}
where~$b_2^+$ denotes the positive part of~$b_2$.
The operator~$\op_\mu$ is self-adjoint and bounded from below
because the form $\opf_\mu$ is densely defined, closed, symmetric
and bounded from below~\cite[Sec.~VII.1.1]{Edmunds-Evans}.
We also verify that
\begin{equation}\label{representation.bis}
  \op_\mu\psi = \ope_\mu \psi \,,
  \qquad
  \psi\in\Dom(\op_\mu) =
  \left\{\psi\in\Dom(\opf_\mu)\,|\,
  \ope_\mu\psi\in\s^2(\Omega) \right\}.
\end{equation}
It is clear that the operators verifying
(\ref{representation}) or~(\ref{representation.bis}) coincide
provided $a\in\s^\infty(\Omega)$ and $b_2=0$.

We recall the definition~(\ref{Rayleigh}) 
of the sequence $\{\gamma_n(\mu)\}_{n\in\Nat^*}$ and prove 
\begin{Theorem}\label{Thm.Asymptotes.App}
Assume~$\Hbasicbis$.
One has
\begin{equation*}
  \forall n\in\Nat^*, \qquad
  \lim_{\mu\to+\infty} \frac{\gamma_n(\mu)}{\mu} 
  = a_\mathrm{min} 
  \,.
\end{equation*}
\end{Theorem}
\begin{proof}
Since 
$
  \op_\mu \geq (\inf\sigma(\op_0) + \mu\,a_\mathrm{min}) \id
$,
we have 
$
  \lim_{\mu\to+\infty} \gamma_n(\mu)/\mu
  \geq a_\mathrm{min} 
$
by the minimax principle.
To obtain the opposite estimate, 
we are inspired by~\cite[proof of Thm.~2.2]{Allegretto-Mingarelli}.
Let $M_a$ denote the maximal operator of multiplication by~$a$, \ie,
$
  M_a \psi := a \psi,  
$ 
$
  \psi \in \Dom(M_a) := 
  \{\psi\in\s^2(\Omega) | \ a\psi\in\s^2(\Omega) \}
$.
The spectrum of~$M_a$ is purely essential 
and equal to the essential range of the generated function~$a$.
In particular, $a_\mathrm{\min}\in\sigma_\mathrm{e}(M_a)$.
By the spectral theorem,
there exists a sequence $\{\psi_j\}_{j\in\Nat^*}\subseteq\Dom(M_a)$
orthonormalized in $\s^2(\Omega)$ such that
$
  \|(M_a-a_\mathrm{min}\id)\psi_j\| \to 0
$
as~$j\to\infty$.
Since $\Smooth_0^\infty(\Omega)$ is dense in $\Dom(M_a)$,
it follows that there is also a sequence 
$\{\varphi_j\}_{j\in\Nat^*}\subseteq\Smooth_0^\infty(\Omega)$
satisfying
$$
  (\varphi_i,\varphi_j)-\delta_{ij} \longrightarrow 0
  \qquad\mbox{and}\qquad
  \big(\varphi_i,(M_a-a_\mathrm{min}\id)\varphi_j\big) \longrightarrow 0
$$ 
as $i,j\to\infty$.
Now, given $N\in\Nat^*$,   
choose $k=k(N)\in\Nat$ sufficiently large so that 
\begin{equation*}
  A(N) - v_\mathrm{min}\,\id \leq N^{-1} \id
  \qquad\mbox{and}\qquad
  B(N) \geq \mbox{$\demi$} \, \id
  \,,
\end{equation*}
where $A(N)$, respectively $B(N)$, is the $N \times N$ symmetric matrix
with entries $(\varphi_{i+k},a\varphi_{j+k})$, 
respectively $(\varphi_{i+k},\varphi_{j+k})$,
for $i,j\in\{1,\dots,N\}$.
In view of~(\ref{Rayleigh}) and since $\Span\{\varphi_{j+k}\}_{j=1}^N$
is an $N$-dimensional subspace of~$\Dom(\opf_\mu)$, 
one has 
$
  \gamma_n(\mu) \leq c_n(\mu,N) 
$
for all $n \in \{1,\dots,N\}$, where $\{c_n(\mu,N)\}_{n=1}^{N}$
are the eigenvalues 
(written in increasing order and repeated according to multiplicity)
of the matrix 
$
  C(\mu,N) \equiv \big(C_{ij}(\mu,N)\big)\mbox{$_{i,j=1}^{N}$}
$ 
defined by
$
  C_{ij}(\mu,N) := 
  \opf_\mu(\varphi_{i+k},\varphi_{j+k})
$.
However, it is easy to see that 
$$
  C(\mu,N)
  \leq 
  \mu \left(a_\mathrm{min}+N^{-1}\right) \, \id
  + d(N) \, \id 
  \,, 
$$ 
where $d(N)$ denotes the maximal eigenvalue of
the matrix 
$$
  \big(
  (\nabla\varphi_{i+k},\nabla\varphi_{j+k})+(\varphi_{i+k},b\varphi_{j+k})
  \big)_{i,j=1}^N
  \,.
$$
Hence, for all $n \in \{1,\dots,N\}$,
$$
  \lim_{\mu\to+\infty} \frac{\gamma_n(\mu)}{\mu}
  \leq a_\mathrm{min}+N^{-1}
$$
with $N$ being arbitrarily large.
\end{proof}
\begin{Remark}
The asymptotic behaviour of Theorem~\ref{Thm.Asymptotes.App}
for both bounded~$\Omega$ and~$a$
was established in~\cite[Corol.~2.6]{Binding-Browne_1987}
(see also \cite[proof of Lemma~2.1]{Binding-Browne_1989}
and \cite{Binding-Browne-Huang}). 
\end{Remark}
%

\section{The semigroup property}\label{App.C0}
%
Unable to find a suitable reference, 
in this appendix we prove that the operator~$A_\alpha$
is the infinitesimal generator of a $C_0$-semigroup on~$\Hilbert$. 
Our strategy is to modify the approach of~\cite[Sec.~XIV.3]{Yosida}
(see also~\cite[Sec.~7.4]{Pazy}),
where the property was shown for~$\Omega=\Real^d$, 
$b\in\s^\infty(\Omega)$ and~$a=0$. 
  
We will need some preliminaries.
Firstly, we prove the following solvability result,
which is well known if~$b$ and~$\Omega$ were bounded.
\begin{Lemma}\label{Lem.solvability}
Assume~$\Hbasic$.
Let~$\eta\in(0,\eta_\alpha)$ with
$
  \eta_\alpha
  := \mbox{$\demi$} \big(1+\alpha|a_\mathrm{min}|+|b_\mathrm{min}|\big)^{-1}
$.
For any $\phi\in\s^2(\Omega)$, 
there is a unique function~$\psi\in\Dom(\op_0)$ satisfying
\begin{equation}\label{solvability}
  \big(1 + \eta\,\alpha a + \eta^2 \op_0\big) \psi = \phi
  \,.
\end{equation}
\end{Lemma}
\begin{proof}
We are inspired by~\cite[Remark~2.3.2]{Ladyzhenskaya}. 
First of all, if~$\eta\in(0,\eta_\alpha)$,
then we can establish the \emph{a priori} bound
\begin{equation}\label{apriori}
  \|\psi\|_{\opf_0} \leq (\sqrt{2}\,\eta)^{-1} \, \|\phi\|
  \,.
\end{equation}
Let~$\{\Omega_n\}_{n\in\Nat}$ be an exhaustion sequence of~$\Omega$
(\ie\/ each~$\Omega_n$ is a bounded open set
such that $\Omega_n\subset\subset\Omega_{n+1}$
and $\bigcup_{n\in\Nat}\Omega_n=\Omega$). 
Let~$\psi_n\in\Hilbert_0^1(\Omega_n)$ 
be the solution of the boundary problem~(\ref{solvability})
restricted to~$\Omega_n$;
the existence and uniqueness of such solutions 
is guaranteed, \eg, by~\cite[Sec.~II.3]{Ladyzhenskaya}.  
We extend~$\psi_n$ to the whole~$\Omega$
by taking it to be zero outside~$\Omega_n$,
and we retain the same notation~$\psi_n$ for this extension. 
All the~$\psi_n$ are elements of~$\Dom(\opf_0)$
and satisfy the \emph{a priori} bound~(\ref{apriori}).
Consequently, $\{\psi_n\}_{n\in\Nat}$ is uniformly bounded in~$\Dom(\opf_0)$
and there exists a subsequence $\{\psi_{n_k}\}_{k\in\Nat}$
which converges weakly in~$\Dom(\opf_0)$ to some function~$\psi$.
This function is the desired solution.
Indeed, the functions~$\psi_n$ satisfy the integral identities
\begin{equation}\label{weak.identity}
  (\varphi,\psi_n) 
  + \eta \, (\varphi,\alpha a\psi_n) 
  + \eta^2 (\varphi,\op_0\psi_n)
  = (\varphi,\phi)
\end{equation}
for all $\varphi\in\Smooth_0^\infty(\Omega_n)$.
In~(\ref{weak.identity}),
we fix some $\varphi\in\Smooth_0^\infty(\Omega_n)$,
and then take the limit in the subsequence~$\{n_k\}_{k\in\Nat}$ 
chosen above, assuming that $n_k \geq n$. 
In the limit we obtain~(\ref{weak.identity})
for~$\psi$ with the function~$\varphi$ that we took.
Since these~$\varphi$ form a dense set in~$\Dom(\opf_0)$,
$\psi$ satisfies the identity~(\ref{weak.identity}) 
for all~$\varphi\in\Dom(\opf_0)$,
hence it is a generalized solution of the problem~(\ref{solvability}).
By the uniqueness, which is guaranteed by the bound~(\ref{apriori}),
the whole sequence~$\{\psi_n\}_{n\in\Nat}$
will converge in this way to the solution~$\psi\in\Dom(\opf_0)$.
The property $\psi\in\Dom(\op_0)$ follows from~(\ref{solvability})
and the facts that~$\phi\in\s^2(\Omega)$ and~$a$ is bounded.  
\end{proof}

Now we can prove
\begin{Lemma}\label{Lem.solution}
Assume~$\Hbasic$.
Let 
$   
  \Phi\equiv\{\phi_1,\phi_2\}
  \in \Smooth_0^\infty(\Omega) \times \Smooth_0^\infty(\Omega)
$
and~$\eta\in(0,\eta_\alpha)$. 
Then the equation
\begin{equation}\label{eq.resolvent}
  \Psi-\eta\,A_\alpha\Psi = \Phi 
\end{equation}
has a unique solution 
$
  \Psi\equiv\{\psi_1,\psi_2\} \in \Dom(A_\alpha)  
$.
Moreover,
$$
  \|\Psi\|_\Hilbert 
  \leq \left(1-\eta/\eta_\alpha\right)^{-1}
  \|\Phi\|_\Hilbert
  \,.
$$
\end{Lemma}
\begin{proof}
Let $\varphi_1,\varphi_2 \in \Dom(\op_0)$ 
be the respective solutions of
\begin{align*}
  \big(1 + \eta\,\alpha a + \eta^2 \op_0\big) \varphi_1 
  &= \phi_1
  \,,
  \\
  \big(1 + \eta\,\alpha a + \eta^2 \op_0\big) \varphi_2 
  &= \alpha a \phi_1 + \phi_2
  \,.
\end{align*}
Set 
\begin{equation*}
  \psi_1 := \varphi_1+\eta\varphi_2
  \,,
  \qquad 
  \psi_2 := -(\alpha a+\eta\, \op_0)\varphi_1+\varphi_2
  \,.
\end{equation*}
It is easy to check that 
$
  \Psi := \{\psi_1,\psi_2\}
  \in \Dom(\op_0)\times\s^2(\Omega)
$ 
and it satisfies the equation~(\ref{eq.resolvent}), \ie,
\begin{equation*}
  \psi_1-\eta\,\psi_2 = \phi_1
  \,,
  \qquad
  \eta \, \op_0 \psi_1 + \big(1 + \eta\,\alpha a\big) \psi_2 = \phi_2
  \,.
\end{equation*}
Actually, it follows from the first of the last two equations
that $\psi_2\in\Dom(\opf_0)$, so~$\Psi\in\Dom(A_\alpha)$ 
and it is the desired solution. 
Moreover, we have
\begin{align*}
  \|\Phi\|_\Hilbert^2  
  \ \equiv \ & \|\phi_1\|_{\opf_0}^2 + \|\phi_2\|^2
  \\
  \ = \ & \|\psi_1-\eta\,\psi_2\|_{\opf_0}^2 
  + \|\eta \, \op_0 \psi_1 + \big(1 + \eta\,\alpha a\big) \psi_2\|^2 
  \\
  \ \geq \ & \|\psi_1\|_{\opf_0}^2 + \|\psi_2\|^2
  + 2\eta \, (\psi_2,\alpha a\psi_2)
  \\
  \ &  
  -2\eta \re\!\big(\psi_2,[-\Delta+b-b_\mathrm{min}+1]\psi_1\big)
  + 2\eta \re(\psi_2,\op_0\psi_1)  
  \\
  \ \geq \ & \|\psi_1\|_{\opf_0}^2 + \|\psi_2\|^2
  + 2\eta \, \alpha a_\mathrm{min} \, \|\psi_2\|^2
  - \eta \left( 
  \|\psi_2\|^2 + |1-b_\mathrm{min}| \, \|\psi_1\|^2 
  \right)
  \\
  \ \geq \ & \left(1-\eta\,\big[
  2+2\alpha|a_\mathrm{min}|+|b_\mathrm{min}|  
  \big]\right)
  \|\Psi\|_\Hilbert^2
  \\
  \ \geq \ & \left(1-\eta/\eta_\alpha\right)
  \|\Psi\|_\Hilbert^2
\end{align*}
and $(1-|x|) \geq (1-|x|)^2$ for~$|x| \leq 1$.
\end{proof}

Lemma~\ref{Lem.solution} shows that the range 
of the operator $\id-\eta A_\alpha$ contains
$\Smooth_0^\infty(\Omega) \times \Smooth_0^\infty(\Omega)$
for all $\eta\in(0,\eta_\alpha)$.
Since~$A_\alpha$ is closed (\cf~Proposition~\ref{Prop.closed}), 
the range of~$\id-\eta A_\alpha$ is all of~$\Hilbert$
and we have 
\begin{Corollary}\label{Corol.solution} 
Assume~$\Hbasic$.
For every $\eta\in(0,\eta_\alpha)$
and $\Phi\in\Hilbert$ the equation~\emph{(\ref{eq.resolvent})} 
has a unique solution 
$
  \Psi \in \Dom(A_\alpha)  
$
and
$
  \|\Psi\|_\Hilbert 
  \leq \left(1-\eta/\eta_\alpha\right)^{-1}
  \|\Phi\|_\Hilbert
$.
\end{Corollary}

Now we are ready to prove
\begin{Theorem}\label{Thm.C0}
Assume~$\Hbasic$.
The operator~$A_\alpha$ is the infinitesimal generator 
of a $C_0$-semigroup $U_\alpha$ on~$\Hilbert$, 
satisfying
\begin{equation}\label{C0.bound}
  \forall t\in(0,\infty), \qquad
  \|U_\alpha(t)\| \leq e^{\omega_\alpha t}
\end{equation}
with
$
  \omega_\alpha :=  
  2 \big(1+\alpha|a_\mathrm{min}|+|b_\mathrm{min}|\big)
$.  
\end{Theorem}
\begin{proof}
The domain of~$A_\alpha$ is clearly dense in~$\Hilbert$
(\cf~Proposition~\ref{Prop.closed}).  
>From Corollary~\ref{Corol.solution} it follows that
$(z-A_\alpha)^{-1}$ exists for all $z\in(\omega_\alpha,\infty)$,
$\omega_\alpha = \eta_\alpha^{-1}$,
and satisfies 
$$
  \big\|(z-A_\alpha)^{-1}\big\| 
  \leq (z-\omega_\alpha)^{-1}
  \qquad\mbox{for}\quad
  z > \omega_\alpha
  \,.
$$
The claim then follows by~\cite[Thm.~1.5.3]{Pazy}. 
\end{proof}
\begin{Remark}
Since~$a$ is bounded, the proof could also be done
by showing that~$A_0$ is an infinitesimal
generator of a $C_0$-group~$U_0$
and applying~\cite[Thm~1.1 of Chap.~3]{Pazy}. 
\end{Remark}

As an application of Theorem~\ref{Thm.C0},
one has the following 
\begin{Corollary}\label{Corol.solution.wave} 
Assume~$\Hbasic$.
Let $\phi_1\in\Dom(\op_0)$ and $\phi_2\in\Dom(\opf_0)$.
Then there exists a unique 
$
  t \mapsto \psi(t,\cdot) 
  \in  \Smooth^0\big([0,\infty);\Dom(\op_0)\big)
  \cap \Smooth^1\big([0,\infty);\Dom(\opf_0)\big)
  \cap \Smooth^2\big([0,\infty);\s^2(\Omega)\big)
$  
satisfying the initial value problem~\emph{(\ref{wave.eq})}. 
\end{Corollary}
\begin{proof}
Let~$U_\alpha$ be the semigroup generated 
by~$A_\alpha$ and set
$$
  \big\{\psi_1(t,x),\psi_2(t,x)\big\}
  := U_\alpha(t) \, \big\{\phi_1(x),\phi_2(x)\big\}
  \,.
$$
Then 
$$
  \frac{\partial}{\partial t} \{\psi_1(t,\cdot),\psi_2(t,\cdot)\}
  = A_\alpha\big\{\psi_1(t,\cdot),\psi_2(t,\cdot)\big\}
  = \big\{\psi_2(t,\cdot),-\op_0\psi_1(t,\cdot)-\alpha a \psi_2(t,\cdot)\big\}
$$
and~$\psi_1$ is the desired solution. 
\end{proof}
%

\section*{Acknowledgements}
\addcontentsline{toc}{section}{Acknowledgements} 
This work was partially supported by FCT/POCTI/FEDER, Portugal,
and the AS\,CR project K1010104.

%
%
\providecommand{\bysame}{\leavevmode\hbox to3em{\hrulefill}\thinspace}
\providecommand{\MR}{\relax\ifhmode\unskip\space\fi MR }
\providecommand{\MRhref}[2]{%
  \href{http://www.ams.org/mathscinet-getitem?mr=#1}{#2}
}
\providecommand{\href}[2]{#2}

\end{document}